\documentclass[12pt,epsf]{article}
\def\reff#1{(\ref{#1})}
\evensidemargin \oddsidemargin

\def\E{{\mathbb E}}
\def\P{{\mathbb P}}
\def\R{{\mathbb R}}

\def\N{{\mathbb N}}

\newcommand{\bx}{\mbox {\protect \boldmath $x$}}
\usepackage{amssymb}
\usepackage{graphicx}
\begin{document}

\title{Detection of spatial pattern through independence of thinned processes} 
\author{
Renato M. Assun\c{c}\~{a}o \\
Departamento de Estat\'{\i}stica \\
UFMG - Universidade Federal de Minas Gerais \\
30161-970 Belo Horizonte, MG, Brazil \\  
assuncao@est.ufmg.br\\
\\
Pablo A. Ferrari \\
Departamento de Estat\'{\i}stica \\
USP - Universidade de S\~{a}o Paulo \\
Caixa Postal 66281 \\
05315-970 S\~{a}o Paulo, SP, Brazil\\
pablo@ime.usp.br} 
\date{} \maketitle

\begin{abstract}
  Let $N$, $N_1$ and $N_2$ be point processes such that $N_1$ is obtained from
  $N$ by homogeneous independent thinning and $N_2 = N- N_1$. We give a new 
  elementary proof that $N_1$ and $N_2$ are independent if and only if $N$ is a 
  Poisson point process. We present some applications of this result to test if a
  homogeneous point process is a Poisson point process.

\end{abstract}

\section{Introduction}

In spatial statistics, it is common to consider simultaneously two
spatial point patterns. A common applied setting is that where the 
researcher considers a pattern composed by the location of disease 
cases in a planar region, and another set of locations labeled as control 
individuals. Generally, the location of case and control individuals
are their residences. 
The attention in the first situation is concentrated on the comparison of 
the marginal distributions of the two processes.
Usually, the interest is to decide if the disease
cases have some degree of spatial clustering with respect to the controls'
pattern (Diggle, 1993; Kelsall and Diggle, 1995), specially around putative 
sources of increased risk (Diggle, 1990; Diggle and Rowlingson, 1994). 
If cases and controls exhibit the same spatial pattern, it 
makes sense to consider the null hypothesis that cases and controls 
are independent random samples from the same population at risk. 
This hypothesis implies that, conditionally on the observed locations of 
cases and controls, the events are labeled by the random outcome of 
fliping a coin of constant probability $p$, where $p$ reflects the 
relative sizes of the cases and controls samples. It is usual to carry out
the test conditioned on the observed number of cases of controls. 

Another common situation in spatial statistics is when the interest concentrates 
on testing the independence of two point patterns and therefore attention is directed
to the joint distribution of the processes. For example, the 
researcher could be studying two species of plants in the same region. 
From theoretical reasons or empirical knowledge, the species could be known 
to have quite different spatial configurations. Therefore, there would be no
interest in testing if they arise by randomly labeling an original process. 
In this situation, it is more usual to test either they are independent point 
processes or, alternatively, if there is interaction between the two processes
(Lotwick and Silverman, 1982;  Wiegand {\em et al.}, 2000).
If the independence hypothesis holds, the expected number of 
individuals from one species in a disc centered at $\bx=(x_1,x_2)$ is 
independent of the presence in $\bx$ of an individual from the other species.

Hence, the two hypothesis are considered in very different situations 
and they imply different consequences to the observed point patterns. 
However, these two hypothesis are not exclusive. Suppose that $N$ is a 
Poisson process with intensity function $\lambda(\bx)$ and that
$N = N_1 + N_2$ where $N_1$ is a thinning of $N$ obtained through the
function $p(\bx)=p$, a constant independent of $N$. That is, $N_1$ is 
a random labeling of the $N$ events. 
Then, it is well known that $N_1$ and $N_2$ are
independent Poisson processes with intensities $p \, \lambda(\bx)$ and
$(1-p) \, \lambda(\bx)$, respectively (Cressie, 1991, page 690).

This result raises the question of the converse statement. Consider a 
point process $N_1$ arising as a random thinning of a point process $N$ 
and let $N_2$ be the complementary point process such that $N=N_1+N_2$.
If $N_1$ and $N_2$ are independent point processes, is it true 
that $N$ is a Poisson point process ? The answer is positive and,
since the Poisson process is the only point process with this property, 
this result gives a characterization of this process. 

This characterization result is not well known among spatial statisticians but
it is not new. Srivastava (1971) proved this characterization for the particular
case of stationary point processes evolving in time.  He provided a short
proof using two previous results: a Poisson distribution characterization by
Moran (1952) and the characterization of a Poisson process by the Poisson
distribution on compact sets by R\'{e}nyi (1967). Fichtner (1975) extended
Srivastava's (1971) characterization theorem for non-stationary point
processes occuring in $\R^d$. 

In this paper, we present a new proof of this characterization of Poisson
processes, possibly non-stationary. We believe our proof is simpler than
Fichtner's. It only uses Moran's theorem and well known point processes
results. We also present a new and elementary proof of Moran's theorem without
using characteristic functions (Lemma 1 below). Since this important
characterization theorem is absent even from major point processes references,
such as Daley and Vere-Jones (1988), we think it will be useful for spatial
statisticians to present it here.

Based on this characterization theorem, we present a new approach to 
test for spatial pattern in an observed point process. Although the 
theorem characterizes also inhomogeneous Poisson processes, in this paper 
we concentrate on the detection of homogeneous Poisson process. We present two
different tests, one based on the bivariate $K$ function, and another based on 
empty space methods. 

We give the definitions and set the notation in Section 2 where
we also prove our main result concerning the characterization of the 
Poisson process through the independence of the processes formed 
by randomly labeling an initial process. In Section 3, we discuss 
some implications for statistical inference about homogeneous
Poisson point processes and we finish with discussion and conclusions
in Section 4.
 
\section{The characterization of a Poisson point processes}

Let $N$ be a point process in $\R^d$ with locally
finite intensity $\nu$: $\nu(K)<\infty$ for each compact set $K$. Let
$N_1$ and $N_1$ independent thinnings of $N$ with acceptance value $p$ and
$1-p$ respectively, with $p\in[0,1]$. These processes are characterized by
\begin{equation}
  \label{119}
  \P(N_1(K) = i, N_2(K) = j) = \P(N(K) = i+j) {i+j\choose i} p^i(1-p)^j
\end{equation}
for any compact set $K$.

\paragraph{Theorem 1} {\sl
The process $N$ is a Poisson process if and only if $N_1$ and $N_2$ are
independent.}

The proof of the theorem is based on an elementary lemma about Poisson
random variables. Let $Z$ be a random variable with values in $\N$ and finite mean
$\lambda>0$. Let $(U_i\,:\,i\in \N)$ be a sequence of independent 
random variables and independent of $Z$ with Bernoulli distribution:
\begin{equation}
\label{100}
  \P(U_i=1) = 1-\P(U_i =0) = p
\end{equation}
where $p\in [0,1]$ is a parameter.

Let $X$ and $Y$ be thinnings of $Z$ using $U_i$:
\begin{equation}
  \label{101}
  X := \sum_{i=1}^Z U_i\;;\quad Y := \sum_{i=1}^Z (1-U_i)
\end{equation}
By Wald identity
\begin{equation}
  \label{113}
  \E X = \lambda p\;;\qquad \E Y = \lambda (1-p)\,.
\end{equation}
Write $r_k= \P(Z=k)$, $p_k= \P(X=k)$ and $q_k= \P(Y=k)$. Then, by
definition:
\begin{eqnarray}
  \label{106}
 &&p_i =  \sum_{n\ge i} r_n {n\choose i} p^i(1-p)^{n-i}\;
;\qquad q_j=\sum_{n\ge  j} r_n {n\choose j} p^{n-j}(1-p)^j
\end{eqnarray}
and
\begin{equation}
  \label{111}
 P(X=i, Y=j) = r_{i+j} {{i+j}\choose{i}} p^i (1-p)^{j}\;,\quad i,j\ge 0\,.
\label{eqxy1}
\end{equation}

\paragraph{Lemma 1}
  {\sl The variable $Z$ has Poisson distribution if and only if $X$ and $Y$ are
  independent.}

\noindent{\bf Proof.} The implication ``$Z$ Poisson implies $X$ and $Y$
independent'' is in textbooks (Cressie, 1991, page 690, for instance). To show
the reverse we first establish the strict positivity of all $r_n$.  Since $X$
and $Y$ are independent, then
 \begin{equation}
   \label{7}
   r_n = \sum_{i+j=n} p_iq_j
 \end{equation}
By \reff{106} $r_n> 0$ implies $q_k>0$ and $p_k> 0$ for all $k\le n$. Since
$Z$ is not identically equal to zero, $r_n>0$ for some $n\ge 1$ and hence
$p_0,q_0,p_1,q_1>0$. By \reff{7} $r_0=p_0q_0>0$ and $r_1\ge p_0q_1>0$.  By
induction, fixing $n\ge 1$ and assuming $r_n>0$, we get $r_{n+1} \ge p_n
 q_1 >0$. This shows that $r_n>0$ for all $n\ge 0$.

Using the hypothesis of independence and taking alternatively $i=x$,
$j=y+1$ and then $i=x+1$ and $j=y$ in \reff{eqxy1} we get
\begin{eqnarray}
  \label{112}
  p_x q_{y+1} &=& r_{x+y+1} {{x+y+1}\choose{x}} p^x (1-p)^{y+1} \\
p_{x+1} q_{y} &=& r_{x+y+1}{{x+y+1}\choose{x+1}} p^{x+1} (1-p)^{y}
\end{eqnarray}
from where
\begin{equation}
\label{121}
p_x q_{y+1} (y+1) p = p_{x+1} q_{y} (x+1) (1-p)\,.
\end{equation}
Fixing $x=0$, \reff{121} and the fact that $(q_y)$ is a probability imply
that
$q_y$ must satisfy:
\begin{eqnarray}
q_{y+1} &=& {\frac{q_{y}}{y+1}} \left( {\frac{1-p}{p}}{\frac{p_1}{p_0}}%
\right) ,\quad y\ge0 \\
\sum_{y\ge 0} q_y &= &1
\end{eqnarray}
whose solution is:
\begin{equation}
q_y = {\frac{1}{y!}} \left({\frac{1-p}{p}}{\frac{p_1}{p_0}}\right)^y e^{-{%
\frac{1-p}{p}}{\frac{p_1}{p_0}}},\quad y\ge0
\label{distpoissony}
\end{equation}
Hence $Y$ has Poisson distribution with mean
${\frac{1-p}{p}}{\frac{p_1}{p_0}}$. By \reff{113} this mean also equals
$\lambda(1-p)$. The same argument shows that $X$ is Poisson with mean
${\frac{p}{1-p}}{\frac{q_1}{q_0}}=\lambda p$. Since $X$ and $Y$ are
independent and $Z=X+Y$, $Z$ must be Poisson. $\square$

\paragraph{Proof of Theorem 1.} A point process is completely
determined by the null probabilities $(\P(N(K)=0\,:\, K$ compact$)$ 
(Theorem 7.3.II in page 216 of Daley and Vere-Jones, 1988).
Denoting $Z=N(K)$, $X=N_1(K)$ and $Y=N_2(K)$, we have that $Z,X,Y$ satisfy
the hypothesis of Lemma 1 with $\lambda = \nu(K)$. 
Hence $N$ is Poisson with intensity $\nu$. $\square$

\section{New tests for homogeneous Poisson point processes}

This characterization of the Poisson process suggests a different way
to test if a point process is a stationary Poisson process. Assume $N$
is a stationary process and, using a coin with success probability
$p$, randomly label some of its events with mark $1$, the remaining events 
being marked as $2$. Only the stationary Poisson process has the two marked 
processes independent. Therefore, to test if the randomly labelled processes 
$N_1$ and $N_2$ are independent is equivalent to test the hypothesis that $N$ 
is a stationary Poisson process.

The usual way to test if two stationary processes observed in a
finite sampling window $A$ with area $|A|$ are independent is that proposed
by Lotwick and Silverman (1982) based on conditional Monte Carlo
tests (Ripley, 1977; Besag and Diggle, 1977) when $A$ is a rectangle.
Firstly, a suitable test 
statistic is chosen reflecting a particular alternative hypothesis of interest. 
If no specific alternatives are enviosioned, it is common to consider the 
bivariate Ripley's $K$ function defined for $d > 0$ by
\begin{equation}
K_{12} (d) = \frac{2 \pi}{\lambda_1 \lambda_2} \;
  \int_0^d u \, \lambda_{12}(u) \, du \;
\end{equation}
where $\lambda_i$ is the first-order intensity of process $N_i$ and
$\lambda_{12}(u)$ is the second-order intensity function of processes
$N_1$ and $N_2$. It is clear that $K_{12}(d)=K_{21}(d)$.
{From} the definitions, it follows that, under independence of $N_1$ and 
$N_2$, we have $K_{12}(d) =  \pi d^2$, whatever the marginal distributions of
the two processes.

Let $n_1$ and $n_2$ be the number of events of $N_1$ and $N_2$, respectively,
observed in the rectangular sampling window $A$. Convert the rectangle to a 
torus by identifying the opposite edges of $A$. With this toroidal idea,
no edge correction is necessary in the definition of the estimator of 
$K_{12}(d)$. Define $I_d(u)$ to be $1$ if $u \leq d$, and $0$ otherwise. 
Let $u_{ij}$ be the distance from the $i$th $N_1$-type event located at 
$\bx_{1i}$ to the $j$th $N_2$-type event. 

The test statistic is based on the empirical function $\tilde{K_{12}}(d)$,
first proposed by Hanisch and Stoyan (1979), and defined by
\begin{equation}
\tilde{K_{12}}(d) = (n_1 n_2)^{-1} |A| \sum_{i=1}^{n_1} \sum_{j=1}^{n_2}
     I_d(u_{ij})
\label{K12def}
\end{equation}
The equality $K_{12}(d) = K_{21}(d)$ is also valid for its 
empirical counterparts $\tilde{K_{12}}(d)$ and $\tilde{K_{21}}(d)$
in this case of a toroidal region. 

Keeping the $N_1$ process fixed, randomly shift the observed 
$N_2$ pattern in the torus and recalculate $\tilde{K_{12}}(d)$. 
After many independent shifts, we have the empirical distribution of 
$\tilde{K_{12}}(d)$ under independence of the processes conditioned
on the observed marginal structure. Percentiles from this distribution
for several different values of $d$ can be used to construct acceptance
envelopes for the hypothesis.      

In the procedure we are proposing, the $N_1$ events are chosen out of 
those from the $N$ process independently with probability $p$. To choose
the value of $p$, consider the variance of (\ref{K12def}). If $N_1$ and
$N_2$ are independent Poisson processes then, conditionally on the
values of $n_1$ and $n_2$,
\begin{eqnarray}
Var\left(\tilde{K_{12}}(d)\right) &=& (n_1 n_2)^{-2} |A|^2
 \sum_{i,i^*=1}^{n_1} \sum_{j,j^*=1}^{n_2}
     Cov \left( I_d(u_{ij},I_d(u_{i^* j^*}) \right) \\ \nonumber
 &=& (n_1 n_2)^{-2} |A|^2  \sum_{i=1}^{n_1} \sum_{j=1}^{n_2}
     Var \left( I_d(u_{ij}) \right) \\  \nonumber 
 &=& (n_1 n_2)^{-1} |A|^2 Var \left( I_{[|X-Y| \leq d]} \right)
\label{VarK12}
\end{eqnarray}  
where $X$ and $Y$ are independent random variables uniformly distributed
over $A$ identified with the torus (Silverman, 1978). 
If $n_1+n_2=n$ is fixed, the optimal choice of $n_1$ and $n_2$ in 
the sense of minimizing the variance ($\ref{VarK12}$) is given by 
$n_1=n_2=n/2$. This suggests labeling the processes with $p=0.5$.

Another possible test statistic is based on the avoidance set function 
$P(N(A)=0)$ or ``empty space'' techniques. There are examples of ergodic 
stationary dependent bivariate point processes that are judged independent 
by second-order methods, such as the $\tilde{K_{12}}(d)$ function, 
but with interactions detected by the avoidance function (Lotwick, 1984).
This leads to the consideration of another test. 

Let $G_1(d)$ be the probability that a disc of radius $d$ contains 
no events of the $N_1$ process. Define $G_2(d)$ and $G(d)$ similarly
for the processes $N_2$ and $N=N_1 + N_2$, respectively. 
If $N_1$ and $N_2$ are independent processes we have the 
following identity holding for all 
$d$: $G(d) = G_1(d) \, G_2(d)$. As a consequence, we can use 
the following statistic to investigate the interaction between 
the processes $N_1$ and $N_2$: 
\begin{equation}
T(d) = \log \hat{G}(d) - \log \hat{G_1}(d) - \log \hat{G_2}(d)
\label{Tddef}
\end{equation}
As previously described, a conditional Monte Carlo test is used to 
assess the significance of empirical estimates of $T(d)$.
Lotwick and Silverman (1984) use the Green-Sibson Dirichlet tesselation algorithm 
for computing the function estimates while we prefer to estimate them from $m$ 
randomly distributed sample points in $A$ as decribed in Diggle (1983, 
page 20).

To choose the value of $p$, consider the variance of (\ref{Tddef}). 
Assuming that $N$ is a Poisson process with intensity $\lambda$ and 
$n$ observed events and ignoring boundary effects, we use a standard 
delta method argument (Taylor expansion) to find  
\begin{equation}
Var(T(d)) \approx \frac{1}{n} \left( e^{2 \lambda \pi d^2} + 1 - 
        (e^{p2 \lambda \pi d^2} + e^{(1-p)2 \lambda \pi d^2} \right) \: .
\label{VarTd}
\end{equation}  
It is clear that the variance is zero when $p=0$ or $p=1$. The reason
is that, in this case, $T(d)=0$ because either $N = N_1$ or $N=N_2$. 
Since $\log \hat{G}(d)$ is fixed whatever value of $p$ is 
chosen, a better strategy is to select $p$ to minimize the variance
of the $log(G_1(d) \, G_2(d))$ estimator. Hence,
\begin{equation}
Var \left( \log \hat{G_1}(d) - \log \hat{G_2}(d) \right)
     \;  \approx\;   \frac{1}{n} \left( 
        (e^{p2 \lambda \pi d^2} + e^{(1-p)2 \lambda \pi d^2} - 2 \right) \: .
\label{VarTd2}
\end{equation}
which is minimized when $p=0.5$, giving a minimum of 
$2/n(\exp(\pi \lambda d^2) -1) > 0$, if $d > 0$.
As we found previously, this new result also suggests 
to label the processes using $p=0.5$.

\subsection*{Example}

%\begin{figure}
%  \begin{center}
%\fbox{\includegraphics[width=150mm]{fig1.eps}}    
%  \end{center}
%\caption{The first, second and third 
%columns of plots refer to the pine saplings, redwood seedlings, and 
%cell centers datasets, respectively. The first, second and third rwos
%of plots refer to the point patterns, the bivariate $K$ functions and 
%the $T(d)$ statistic, respectively. The dashed lines are approximate 
%95\% confidence bands.}
%\end{figure}

We illustrate the techniques described with some real data:  
the locations of 62 
redwood seedlings in a square of 23 meters, the locations of
42 biological cell centers in a unit square, and 
the locations of 65 Japanese black 
pine saplings in a square of side 5.7 meters. All the data are as 
reported by Diggle (1983) from the references therein. 

\fbox{Figure 1 around here}

Figure 1 shows the results of our two tests, based on the $K_{12}(d)$
function and in the empty space function. The first, second and third 
columns of plots refer to the redwood seedlings, cell centers, 
and pine saplings, respectively. The first row of plots shows the three point 
patterns. We used $p=0.5$ to generate the thinned processes showed as 
circles and crosses in Figure 1. The second and third rows of plots 
refer to the $\hat{K_{12}}(r)$ and $\hat{T}(d)$ tests, respectively.  

Several tests have been used previously in these datasets and usually they
accept the hypothesis of a homogeneous Poisson process for the pine saplings,
and reject this hypothesis for the clustered redwood seedlings pattern and the
regularly spaced cell centers. Our tests find these same results as can be
seen by the behavior of the observed $\hat{K_{12}}(r)$ and $\hat{T}(d)$
functions with respect to the 95\% confidence envelopes. Both test functions
lie outside the envelopes for the first and second datasets and inside the
envelope for the third dataset.

\section{Discussion and conclusions}

A fundamental property characterizing the Poisson process is the
independence of counts on disjoint areas. 
The characterization theorem presented in this paper
suggests that independence of random partitions of events
in the same area is also capable of characterizing the Poisson process.
This is another justification for the usual labeling of a 
homogeneous Poisson process as {\em complete spatial
randomness} (Diggle, 1983). 

A result related to this theorem is Raikov's theorem (see Daley and Vere-Jones,
1988, page 31) which shows that if $Z$ is a Poisson random variable 
expressible as a sum $Z=X+Y$ of independent nondegenerate, nonnegative 
random variables then $X$ and $Y$ are Poisson random variables. The 
present characterization theorem 
drops the hypothesis that $Z$ has a Poisson distribution and shows
that this is a consequence of the independence of $X$ and $Y$ if they are 
obtained through the thinning of $Z$. 
 
We used the characterization result to propose two tests, based on
empty space and second-order methods, for
the hypothesis that a point process is a stationary Poisson process. 
It has not been considered in this paper their relative power 
in detecting departures from the null hypothesis
of a Poisson process and the kind of departure detected by the two techniques.
Likewise, we have not considered the relative merits of other techniques such as 
that based on the $K$ function of the single type $N$ process.

Although the characterization 
result is also valid for non-homogeneous Poisson point processes, it 
is not clear how this could be used to set up a hypothesis test in this 
case. The main problem in the non-homogeneous situation is the 
dependence of the thinned processes on the unknown first-order
intensity function of $N$. Similar problems have made difficult 
the estimation of second moment functions for stationary Cox process
(Chetwynd and Diggle, 1998).  

The examples in Section 3 demonstrate that tests for interaction between two
complementary processes obtained through the random thinning of a stationary point
process $N$ can be an alternative to test the hypothesis that $N$ is a stationary 
Poisson process. Both tests, that based on the $K_{12}$ function and that based 
on empty space methods, lead to the same conclusions in the examples considered
in this paper. These conclusions are the same reached using other tests previously 
proposed in the literature. 

\section*{Acknowledgments}

This research was carried out while the first author was in a sabatical leave
visiting both, Funda\c{c}\~{a}o Oswaldo Cruz, FIOCRUZ, and Escola Nacional 
de Ci\^{e}ncias Estat\'{\i}sticas, ENCE-IBGE, Rio de Janeiro, Brazil, whose 
hospitality is gratefully acknowledged. The visit was supported by  
FORD Foundation grant 990-1161. 
The research received also support from CNPq grant 465928/2000-5. 
%The authors express thanks to FULANO for helpful comments on an earlier 
%version of the paper, to BELTRANO for assitance in computational
%strategies, and to an anonymous referee for helpful comments.

\section{References}

\begin{itemize}

\item Besag, J. and Diggle, P. J. (1977) Simple Monte Carlo tests for 
spatial patterns. {\em Applied Statistics}, {\bf 26}, 327-333.

\item Chetwynd, A. G. and Diggle, P. J. (1998) On estimating the reduced 
second moment measure of a stationary point process. {\em Australian 
and New Zealand Journal of Statistics}, {\bf 40}, 11-15.

\item Cressie, N. (1991) {\em Statistics for spatial data}. New York: 
John Wiley \& Sons.

\item Daley, D.J., and Vere-Jones, D. (1988) {\em An introduction to the 
theory of point processes}. New York: Springer-Verlag.

\item Diggle, P. J. (1983) {\em Statistical Analysis of Spatial Point 
Patterns}. London: Academic Press. 

\item Diggle, P.J. (1990). A point process modelling approach to raised 
incidence of a rare 
phenomenon in the vicinity of a pre­specified point. {\em Journal of the 
Royal Statistical Society A}, {\bf 153}, 349-3­62. 

\item Diggle, P. J. (1993) Point process modelling in environmental epidemiology,
in Barnett, V. and Turkman, K.F. (eds.) {\em Statistics for the Environment}.
Chichester: John Wiley. 

\item Diggle, P.J. and Rowlingson, B. S. (1994). A conditional approach to 
point process modelling of raised incidence. {\em Journal of the 
Royal Statistical Society A}, {\bf 157}, 433-440. 

\item Fichtner, V. K.-H. (1975) Charakterisierung Poissonscher 
zuf\"{a}lliger Punkfolgen und infinitesemale Verd\"{u}nnungsschemata.
{\em Mathematische Nachrichten}, {\bf 193}, 93-104. 

\item Hanisch, K. H. and Stoyan, D. (1979) Formulas for second-order analysis 
of marked point processes. {\em Mathematische Operationsforschung und Statistik,
Series Statistics}, {\bf 10}, 555-560. 

\item Kelsall, J. and Diggle, P.J. (1995). Kernel estimation of relative 
risk. {\em Bernoulli}, {\bf 1}, 3-16. 

\item Lotwick, H. W. (1984) Some models for multype spatial point processes,
with remarks on analysing multitype patterns. {em Journal of Applied 
Probability}, {\bf 21}, 575-582.

\item Lotwick, H. W. and Silverman, B. W. (1982) Methods for analysing 
spatial processes of several types of points. {\em Journal of the Royal 
Statistical Society} B, {\bf 44}, 406-413.

\item Moran, P. A. P. (1952) A characterization of the Poisson distribution.
{\em Proceedings of the Cambridge Philosophical Society}, {\bf 48}, 206-207.

\item R\'{e}nyi, A. (1967) Remarks on the Poisson process. {\em Symposium
on Probability Methods in Analysis}. Berlin: Springer-Verlag, 280-286.

\item Ripley, B. D. (1977) Modelling spatial patterns (with discussion). 
{\em Journal of the Royal Statistical Society} B, {\bf 39}, 172-212.

\item Silverman, B. W. (1978) Distances on circles, toruses and spheres. 
{\em Journal of Applied Probability}, {\bf 15}, 136-143.

\item Srivastava, R. C. (1971) On a characterization of the Poisson
process. {\em Journal of Applied Probability}, {\bf 8}, 615-616.  

\item Wiegand, K., Jeltsch, F. and Ward, D. (2000) Do spatial effects 
play a role in the spatial distribution of desert-dwelling Acacia raddiana?
{\em Journal of Vegetation Science}, {\bf 11}, 473-484.

\end{itemize}

\newpage 

\begin{figure}
  \begin{center}
\fbox{\includegraphics[width=150mm]{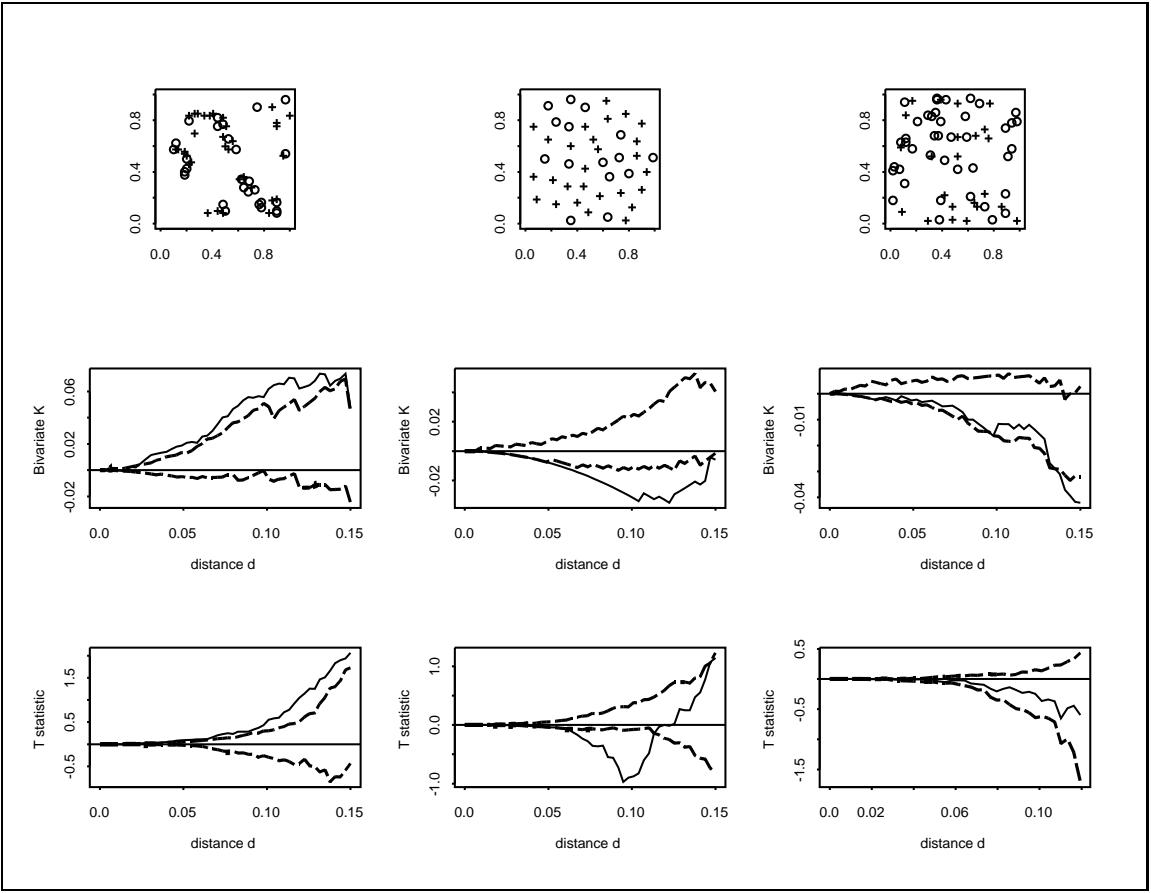}}    
  \end{center}
\caption{The first, second and third 
columns of plots refer to the redwood seedlings, cell centers,
and pine saplings datasets, respectively. The first, second and third 
rows of plots refer to the point patterns, the bivariate $K$ functions 
and the $T(d)$ statistic, respectively. The dashed lines are 
approximate 95\% confidence bands.}
\end{figure}

\end{document}